\documentclass[12pt,bezier]{article}
\usepackage{amsmath,amssymb,amsfonts,euscript,graphicx}

\newtheorem{preproof}{{\bf Proof.}}

\newenvironment{proof}[1]{\begin{preproof}{\rm
               #1}\hfill{$\Box$}}{\end{preproof}}

\title{\large \bf On Maximal and Minimal Linear Matching Property
\thanks{2010 MSC: 12F05. }
\thanks
{{\it Key Words}: Linear matching property, Algebraic number field, Field extension, Maximal linear matching property, Minimal linear matching property}}

\author{{\normalsize{\sc M. Aliabadi${}^{\mathsf{a,1}}$}, {\normalsize{\sc M. R. Darafsheh${}^{\mathsf{b,2}}$}}}\\
\vspace{3mm}
\\{\footnotesize{${}^{\mathsf{a}}$\it Department of Mathematical Sciences, Sharif University of Technology, Tehran, Iran, }}
{\footnotesize{}}\\{\footnotesize{${}^{\mathsf{b}}$\it School of Mathematics, Statistics and Computer Science, Colledge of Science, University of Tehran, Tehran, Iran.}}\\
{\footnotesize{}}\\
{\footnotesize{$^{1}{E-mail~ address: mohsenmath88@gmail.com}$\hspace*{1cm}
$^{2}{E-mail~ address: darafsheh@ut.ac.ir}$
}}}

\date{}
\begin{document}

\maketitle

\begin{abstract}
{\small \noindent 
The matching basis in field extentions is introduced by S. Eliahou and C. Lecouvey in \cite{2}. In this paper we define the minimal and maximal linear matching property for field extensions and prove that if $K$ is not algebraically closed, then $K$ has minimal linear matching property. In this paper we will prove that algebraic number fields have maximal linear matching property. We also give a shorter proof of a result established in \cite{6} on the fundamental theorem of algebra.}
\end{abstract}

\vspace{9mm}
\section{Introduction}
Throughout this paper we will consider a field extension $K\subset L$ where $K$ is commutative and central in $L$. Let $G$ be an additive group and $A,B\subset G$ be nonempty finite subsets of $G$. A {\it{matching}} from $A$ to $B$ is a map $\phi: A \rightarrow B$ which is bijective and satisfies the condition
\[
a + \phi(a) \not\in A
\]
for all $a \in A$. This notion was introduced in \cite{3} by Fan and Losonczy, who used matchings in $\mathbb Z^n$ as a tool for studying an old problem of Wakeford concerning canonical forms for symmetric tensors \cite{7}. 
Eliahou and Lecouvey extended this notion to subspaces in a field extension, here we will introduce a notion from \cite{2}.\\ Let $K\subset L$ be a field extension and $A, B$ be $n$-dimensional $K$-subspaces of $L$. Let $\mathcal A=\{a_1,\ldots,a_n\}$, $\mathcal B=\{b_1,\ldots,b_n\}$ be basis of $A$ and $B$ respectively. It is said that $\mathcal A$ is {\it{matched}} to $\mathcal B$ if
\[
a_i b\in A~~\Rightarrow~~b\in\langle b_1,\ldots,\hat{b_i},\ldots,b_n\rangle
\]
for all $b\in B$ and $i=1,\ldots,n$, where $\langle b_1,\ldots,\hat{b_i},\ldots,b_n\rangle$ is the hyperplane of $B$ spanned by the set $\mathcal B\setminus\{b_i\}$. Also it is said that $A$ is {\it matched} to $B$ if every basis of $A$ can be matched to a basis of B.
\\
It is said that $L$ has the {\it linear matching property} from $K$ if, for every $n\geq1$
and every $n$-dimensional $K$-subspaces $A$ and $B$ of $L$ with $1\not\in B$, the subspace $A$ is matched to $B$. By this we mean linear matching property for $K$-subspaces.
\\
As we mentioned, the above notion was introduce by Eliahou and Lecouvey in \cite{2}, where they proved that if $K\subset L$ is a field extension and $[L: K]$ is prime, then $L$ has linear matching property (see Theorem 5.3 in \cite{2}). We extend this property to the family of field extensions and introduce the notions of minimal and maximal linear matching properties.
\section{Definitions and the main results}
\textbf{Definition 2.1.} Let $K$ be a field. We say $K$ has {\it minimal linear matching property} if there exists a finite field extension $L$  of $K$, such that $L$ has linear matching property from $K$.
\\
\textbf{Definition 2.2.} Let $K$ be a field. We say $K$ has {\it maximal linear matching property} if for any positive integer $n$, there exists a field extension $L_n$ of $K$, such that $[L_n: K]$= $n$ and $L_n$ has linear matching property from $K$.
\\
We shall prove the following results in section 5.
\\
\textbf{Theorem 2.3.} Let $K$ be a field which is not algebraically closed, then $K$ has the minimal linear matching property.
\\
\textbf{Theorem 2.4.} Algebraic number fields have the maximal linear matching property.
\\
\textbf{Theorem 2.5.} Suppose that $K$ is a field and has the maximal linear matching property, then $K$ is infinite.
\\
To prove our main results, we will use Theorem 3.1 which can be regarded as an improvement of the foundamental theorem of algebra. In \cite{6}, Shipman gives an algebraic proof of the foundamental theorem of algebra in special cases, but here we present a different proof which is independent Shipman's proof.
\section{An improvment of the fundamental theorem of algebra}
\textbf{Theorem 3.1.} Let $K$ be a field such that every polynomial of prime degree in $K[x]$ has a root in $K$, then $K$ is algebraically closed.
\\
\begin{proof}
{
 First, we claim there exists a prime $p$ such that for any non-linear irreducible polynomial $f(x)\in K[x]$, $p$ divides the degree of $f(x)$. Suppose that this claim is false, and $p_1,\ldots,p_n$ are prime divisors of the degree of $f(x)$, then there exists $g_i\in K[x]$ such that $p_i\not|\deg g_i(x)$ and $g_i(x)$ is an irreducible polynomials in $K[x]$, where $1\leq i\leq n$.
\\
Now set $F(x):=f^{k_0}(x)g_1^{k_1}(x)\cdots g_n^{k_n}(x)$ where $k_0,k_1,\ldots, k_n$ are non-negative integers. It is clear that $\gcd(\deg f(x),\deg g_1(x),\ldots,\deg g_n(x))=1$ and $\deg F=k_0\deg f+k_1\deg g_1+\cdots+k_n\deg g_n$. By Dirichlet's Theorem on primes, since the $k_i$'s are non-negative integers, we can choose $k_0,\ldots,k_n$ such that $\deg F$ becomes a prime number. So $F(x)$ has a root in $K$ and this is a contradiction. Therefore there exists a prime $p$  such that $p$ divide the degree of every irreducible polynomials in $K[x]$. Now if $L$ is a field extension  of $K$ of degree $p$ and $\alpha\in L\setminus K$, then $L=K(\alpha)$ and if $f(x)\in K[x]$ is the minimal polynomial of $\alpha$, then $\deg f(x)=p$ and $f(x)$ has a root in $K$ and this is a contradiction, hence $K$ has no field extension of degree $p$. Let $L$ be a Galois extension of $K$ with $[L: K]=p^r\cdot m$ where $r, m\in\mathbb N$, $(m,p)=1$. By Galois fundamental theorem and Cauchy theorem, there is an intermediate field $L'$, $K\subset L'\subset L$ such that $[L: L']=p^r$, then $[L': K]=m$. If $m>1$ we can choose $\alpha\in L'\setminus K$, and assume $f(x)$ is the minimal polynomial of $\alpha$ over $K$, then $\deg f(x)|m$, also $f(x)$ is irreducible, then $p|\deg f(x)$, so $p|m$, a contradiction. Hence $m=1$ and $[L: K]=p^r$, again by using Galois fundamental theorem and Cauchy theorem there exists an intermediate field $L'$, $K\subseteq L'\subset L$ such that $[L: L']=p^{r-1}$, then $[L': K]=p$, but since we proved that $K$ has no field extension of degree $p$, this is a contradiction. Thus $K$ has no Galois extension and it is algebraically closed.
}
\end{proof}
\textbf{Corollary 3.2}  Let $K$ be a field such that every polynomial of prime degree in $K[x]$ is reducible on $K$. Then $K$ is algebraically closed.\\
\section{Preliminary results about field extensions and linear matching property}
We use the following result from \cite{4}.
\\
\textbf{Theorem 4.1.} Let $L$ be a finite field of characteristic $p>0$ where $\mathbb Z_p$ is embedded in $L$ and $[L:\mathbb Z_p]=n$. Then for any divisor $m$ of $n$, $L$ has a subfield with $p^m$ elements.
\\
We also use the following result from \cite{5} which is about field extensions with no proper intermediate subfield.
\\
\textbf{Theorem 4.2.} If $K$ is an algebraic number field, then for every positive integer $n$ there exist infinitely many field extensions of $K$ with degree $n$ having no proper subfields over $K$.
\\
The following theorem was proved in \cite{2}, see also \cite{1}.
\\
\textbf{Theorem 4.3.} Let $K\subset L$ be a field extension. Then $L$ has linear matching \linebreak{property} if and only if $K\subset L$ has no proper intermediate subfield with finite degree over $K$.
\\
Now we are ready to prove the main results.
\section{Proof of main results}
\textbf{Proof of Theorem 2.3}
\\
\begin{proof}
{
By Corollary 3.2 there exists an irreducible polynomial $f(x)$ of prime degree in $K[x]$. Now if $L$ is the splitting field of $f(x)$ over $K$, then $[L: K]$ is prime and by Theorem 4.3 $L$ has the linear matching property from $K$, so $K$ has the minimal linear matching property.
}
\end{proof}
\textbf{Proof of Theorem 2.4}
\\
\begin{proof}
{
Let $K$ be an algebraic number field. Then by theorem 4.2 for any positive integer $n$, there exists an extension $L_n$ of $K$ with $[L_n: K]=n$ and this field extension has no proper intermediate subfield, then by Theorem 4.3, $L_n$ has the linear matching property from $K$, so $K$ has the maximal linear matching property.
}
\end{proof}
\textbf{Proof of Theorem 2.5}
\\
\begin{proof}
{
Let $K$ be a finite field with $|K|=p^n$ and $p$ a prime and $n$ a positive integer. Now let $q$ and $m$ be positive integers with $n<q<m$ and $q|m$. If $L$ is an extension of $K$ of degree $m$, then $[L: \mathbb{Z}_p]=mn$ and by Theorem 4.1, $\mathbb{Z}_p \subseteq L$ has an intermediate subfield $K'$ of degree $p^q$. Now since finite fields with the same cardinality are isomorphic, $K'$ is a finite proper intermediate subfield in the extension $K \subset L$ with finite degree over $K$, then by Theorem 4.3, $L$ does not have linear matching property from $K$, hence $K$ does not have maximal linear matching property.
}
\end{proof}

{}

\end{document}